\documentstyle[10pt,twocolumn,openbib]{article}

\newcommand{\Real}{\mbox{\boldmath $R$}}

\newcount\eqcount  
\def\eq{\global\advance\eqcount by1
        \eqno{\hbox{(\number\eqcount)}}}

\newcommand{\beqa}{\begin{eqnarray}}
\newcommand{\eeqa}{\end{eqnarray}}
\def\beqarr{\begin{eqnarray}}
\def\eeqarr{\end{eqnarray}}

\newcommand{\beq}{\begin{equation}}
\newcommand{\eeq}{\end{equation}}

\newtheorem{theorem}{Theorem}

\input {epsf.tex}
\begin{document}
\title{\large\bf MATCHING CONTROL LAWS FOR A BALL AND BEAM SYSTEM}
\author{{\bf F. Andreev}\thanks{ Department of Mathematics, 
on leave from Steklov Institute of Mathematics, St.-Petersburg, Russia}, 
{\bf D. Auckly}$^\ddag$, 
{\bf L. Kapitanski}\thanks{ Department of Mathematics }, 
{\bf A. Kelkar}\thanks{ Department of Mechanical and Nuclear Engineering},\\ 
and {\bf W. White} $^\S$ \\{\it Kansas State University, Manhattan, KS 66506}}
\date{\parbox{5.0in}{\small{
Abstract: This note describes a method for generating an 
infinite-dimensional  family 
of nonlinear control laws for underactuated systems.  
For a ball and beam system, the entire family is found explicitly. 
{\em Copyright $\copyright$ 2000 IFAC}}}\\ \vspace*{.05in}
 \parbox{5in}{\begin{flushleft}\small{Keywords: Nonlinear control, mechanical systems}\end{flushleft}}}
\maketitle

\footnotetext[1]{This work was partially supported 
by NSF Grant No. CMS-9813182.}

\thispagestyle{empty}
\vspace*{-0.25in}

\renewcommand{\baselinestretch}{1.0}
       \pagestyle{empty}
\vspace*{-.1in}
\centerline{1. THE MATCHING CONDITION}

\noindent 
This note presents an application of the method developed by 
Auckly, {\it et al.} (2000), to stabilization of a ball and beam system. 
The results are fully described in (Andreev, {\it et al.}, (2000), 
Auckly, Kapitanski (2000), and 
Auckly, {\it et al.} (2000)). An experimental comparison of a linear control
law versus the nonlinear
control laws described here will be given in the full paper, (Andreev, {\it et al.}, (2000)).  

\medskip

\noindent Let $\,Q\,$ denote a configuration space.
Let $g\in\Gamma(T^\ast Q\otimes T^\ast Q)$ be a metric.
Let $c,f:TQ\to TQ$ be fiber-preserving maps. We assume that $c(-X)=-c(X)$. 
Let $V:Q\to{\Real}$. The
differential equation that we consider is
\begin{equation}
 \nabla_{\dot\gamma}\dot\gamma 
+ c(\dot\gamma) +\ grad_\gamma V
   =f(\dot\gamma).   
\label{eq1}
\end{equation}
Let $P\in\Gamma(T^\ast Q\otimes TQ)$ be a $g$-orthogonal projection. 
We consider the situation where a constraint $P(f) = 0$ is imposed.
A system is called underactuated if $P \neq 0$.

\smallskip
\noindent
Several recent papers propose to find control inputs so that 
the closed-loop system (\ref{eq1}) would have a natural 
candidate for a Lyapunov function (Bloch, {\it et al.} (1998), 
Hamberg (1999), and van der Shaft (1986)).
Auckly, {\it et al.} (2000) introduced the following matching condition and 
characterization of matching in terms of linear partial differential equations.
A control input, $f$, satisfies
the matching condition if there are functions $\widehat g$, 
$\widehat c$, and $\widehat V$ so that the closed loop equations
take the form:
\begin{equation}
 \widehat \nabla_{\dot\gamma} \dot\gamma 
+ \widehat c(\dot\gamma) +\ \widehat  {grad}_\gamma
   \widehat V=0. 
\label{eq2}
\end{equation}
The motivation for this method is that 
$
\widehat H = \frac{1}{2}\widehat g(\dot \gamma,\dot \gamma) 
+ \widehat V(\gamma)
$
is a natural candidate for a
Lyapunov function because 
$
d\widehat H /dt=-\widehat g(\widehat c(\dot\gamma),\dot\gamma) $.
A straightforward computation shows that, the matching condition is satisfied 
if and only if
\begin{equation}
P(\nabla_X X-\widehat \nabla_X X)=0,  
\label{eq3}
\end{equation}
\begin{equation}
\begin{array}{l} 
P(grad_\gamma V-\ \widehat  {grad}_\gamma \widehat V)=0, \\  
P(c(X)-\widehat c(X))=0.
\end{array}  
\label{eq4}
\end{equation}
Equation (\ref{eq3}) is a system of non-linear 
first order PDE's for  
$\widehat g$. It is perhaps surprising and pleasing that all of the
solutions to (\ref{eq3}), (\ref{eq4}) may be obtained by 
first solving one first order 
linear system of PDE's and then solving a second set of linear PDE's. 
This is accomplished by 
introducing a new
variable, $\lambda$, by $g(X,Y)=\widehat g(\lambda X,Y)$.

\begin{theorem} 
The metric, $\widehat g$, satisfies (\ref{eq3}) if and only if 
$\lambda$ and $\widehat g$ satisfy
\begin{equation} 
\nabla g\lambda\big|_{\hbox{Im}\ P^{\otimes 2}}=0, \quad
L_{{}_{\lambda PX}}\widehat g = L_{{}_{PX}} g. 
\label{eq5}
\end{equation} 
\end{theorem}

\noindent
In the special case of a system with two degrees of freedom, it is possible
to write out the general solution to this set of differential equations. 
Following Auckly, Kapitanski (2000), express the underactuaded subspace
as the span of a unit length vectorfield, $PX$.  
Choose coordinates $x^1$, $x^2$ so that $PX={\partial\over\partial x^1}$,
and write $\lambda PX=\sigma {\partial\over\partial x^1}
+\mu {\partial\over\partial x^2}$. For the $\lambda$-equation, (\ref{eq5}),
to be consistent 
the following compatibility condition must hold: 
$
\partial ([11,2]\,\mu)/\partial x^2=
\partial ([12,2]\,\mu)/\partial x^1 $.
Starting with this equation and working backwards, all of the equations may
be solved via the method of 
characteristics.
\bigskip

\vspace*{-.05in}
\centerline{2. THE BALL AND BEAM SYSTEM}
\smallskip
\centerline{\epsfxsize=2.29in\epsffile{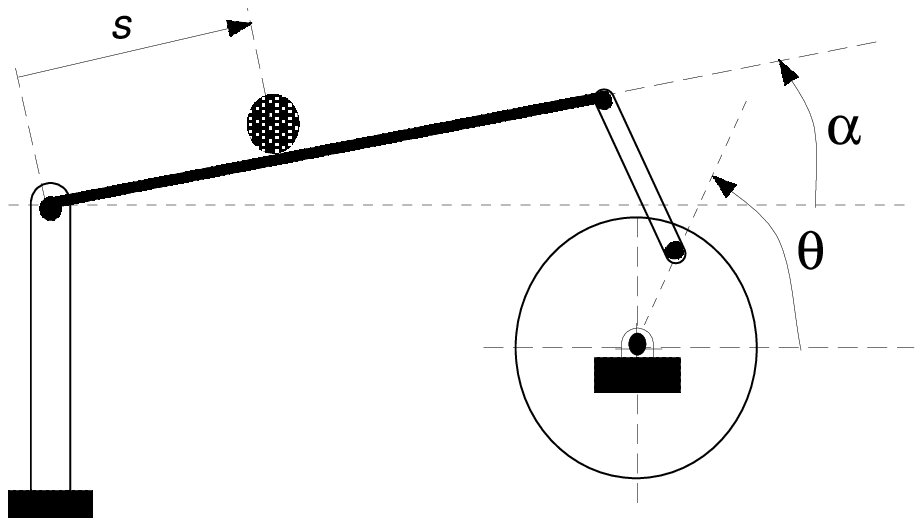}}
\noindent Fig.1. Nonlinear mechanical system.
\smallskip

\noindent As an application of our method consider 
the stabilization problem for the ball and beam system described 
schematically in figure 1. 
One can express $\alpha$ as an explicit function of $\theta$.  
After rescaling, the kinetic energy of the system is given by:
$$
T= \frac12\dot s^{2} 
+  \alpha^\prime
\dot s
\dot \theta
+ \frac12\left (a_4+\left (a_3+\frac 52{s}^{2}\right )
\left(\alpha^\prime\right )^2\right) \dot \theta^2
$$
and 
$
V=a_5\sin(\theta)+(s+a_6)\sin(\alpha)
$,
where the $a_k$ are dimentionless parameters. 
The projection, 
$P=(ds+\alpha^\prime d\theta)\otimes \partial/\partial s$, 
so the control input $u$ is related to $f$ in (1) by 
$f=(u d\theta)^\sharp$.
The resulting equations of motion are
$$
\ddot s +\alpha^\prime \ddot\theta 
+(\alpha^{\prime\prime}-\frac 52 s {\alpha^\prime}^2)\dot\theta^2+
\sin(\alpha)=0 
$$
$$
\begin{array}{l}
\alpha^\prime\ddot s+[a_4+(a_3+\frac 52 s^2) {\alpha^\prime}^2]\ddot\theta
+5{\alpha^\prime}^2s\dot s\dot\theta \qquad\qquad\\ 
\qquad +
(a_3+\frac 52 s^2)\alpha^\prime\alpha^{\prime\prime}\dot\theta^2 
+a_5\cos\theta \\
\qquad \qquad +(a_6+s)\cos(\alpha)\,\alpha^\prime+a_7 \dot\theta=u\,,
\end{array}
$$
where $a_7$ corresponds to inherent dissipation. 
\medskip

\noindent The general solution to the matching equations is
$$
\widehat g_{11}(s,\theta)= 
\psi^2(\alpha)\,(h(y(s,\theta))+ 
 10\int_0^\alpha { d\varphi   \over
\mu_1^{\prime}(\varphi)\psi^2(\varphi)})
$$
$$
\widehat g_{12} ={1\over\mu} (g_{11}-\sigma \widehat g_{11}),
\quad 
\widehat g_{22} ={1\over\mu} (g_{12}-\sigma \widehat g_{12}),
$$
      \renewcommand{\baselinestretch}{.5}
$$
\begin{array}{l}
\widehat V(s,\theta)=w(y)\,+\,
5\,(y+s_0)\,\int_0^\alpha
{\sin(\varphi)\over\mu_1^\prime(\varphi)\psi(\varphi)}\,d\varphi\,
\\
\\
 \qquad\qquad -\,
5\,\int_0^\alpha 
{\sin(\varphi)\over\mu_1^\prime(\varphi)\psi(\varphi)}\,
\int_0^\varphi \psi(\tau)\,d\tau\;d\varphi ,
\end{array}
$$
where 
$
y=\psi(\alpha)s-s_0+\int_0^\alpha \psi(\tau)\,d\tau
$,
$
\psi(\alpha)=
\exp \{-5 \int_0^\alpha {\mu_1(\kappa)\over \mu_1^{\prime}(\kappa)}\,d\kappa\}
$,
$
\mu(s,\theta)={\mu_1^{\prime}(\alpha)\over 5 s\,\alpha^\prime}$,
$
\sigma(s,\theta)=\mu_1(\alpha)-{1\over 5s}\,\mu_1^{\prime}(\alpha)
$
and $\mu_1$, $h$, and $w$ are arbitrary functions. Also, 
$\widehat c^1= -\alpha^\prime \widehat c^2$, where 
$\widehat c^2(s,\theta,\dot s,\dot \theta)$
is an arbitrary function which is odd in $\dot s$ and $\dot\theta$. 
The final nonlinear control law is $u=u_g+u_V+u_c$, where 
$
u_g=g(\nabla_{\dot\gamma}{\dot\gamma}
-\widehat\nabla_{\dot\gamma}{\dot\gamma},{\partial\hfil\over\partial\theta})$, 
$u_V={\partial V\over\partial\theta}-g(\widehat{grad}_\gamma \widehat V,
{\partial\hfil\over\partial\theta})$, and 
$u_c=a_7\,\dot\theta 
- g(\widehat c(\dot\gamma),{\partial\hfil\over\partial\theta})$. 
Using $\widehat H$ as a Lyapunov function, we obtain the following 
conditions that guarantee local asymptotic stability of the equilibrium:
$\det(\widehat g(0))>0$, $\hbox{tr}(\widehat g(0))>0$, $\det(\widehat g 
\widehat c(0))>0$,
$\hbox{tr}(\widehat g \widehat c(0))>0$, $\det(D^2\widehat V(0))>0$, and
$\hbox{tr}(D^2\widehat V(0))>0$,


\medskip

\noindent Another way to check local asymptotic stability is to find 
the poles of the linearized closed-loop system. It is a theorem 
(Andreev, {\it et al.} (2000), Auckly, Kapitanski (2000)) that any linear full state feedback 
control law can be obtained as a linearization of 
some control law in our family.
 
\smallskip
\noindent 
A good stabilizing control law will produce a large basin of 
attraction, send solutions to the equilibrium in a short period 
of time, and will require little control effort. It is, unfortunately, 
not clear how to quantify 
these goals. 

\noindent
We have done some numerical simulation of 
various control laws in our family. We always pick the arbitrary
functions in our nonlinear 
control law in such a way that 
the linearization at the desired equilibrium, 
$u_{lin}=a_8+K_{bp}(s-s_0)+K_{ap}\theta+K_{bd}\dot s+K_{ad}\dot\theta$, 
is exactly the linear 
state feedback control law provided by the manufacturer 
of a commercially 
available system (Apkarian, (1994)).  The numerical and experimental 
response of the system to various initial conditions 
will be recorded in the full version of the paper.

\medskip
\centerline{3. CONCLUSION}

\noindent 
We believe that nonlinear control laws have the potential to achieve better
performance than linear control laws. There are, however, several subtle 
questions which must be resolved before nonlinear control laws may be 
fully exploited in practice. The first question is how to quantify
performance. The second question is how to pick a control law which will
come close to optimizing performance. One interesting idea is to restrict
attention to a class of control laws which generate a closed loop system
of a special form. The hope is then that it will be easier to quantify the
performance of such systems. We have shown that, 
in many situations it is possible to find all control laws which will 
result in a closed loop system  
of the form (\ref{eq2}). 
 
\centerline{REFERENCES}
\baselineskip=11.5pt
\begin{description}
\item Andreev, F., D. Auckly, L. Kapitanski, A. Kelkar, and W. White (2000).
Matching, linear systems, and the ball and beam.
{\it Preprint.}
\vspace*{-.1in}
\item Apkarian, J. (1994). Control System Laboratory, Quanser Consulting, Hamilton, Ontario, Canada L8R 3K8.
\vspace*{-.1in}
\item Auckly, D., L. Kapitanski, and W. White (2000).
Control of nonlinear underactuated systems. 
 To appear in {\it Commun. Pure Appl. Math.}
\vspace*{-.1in}
\item Auckly, D., L. Kapitanski (2000).
Mathematical Problems in the Control of Underactuated Systems.
{\it Preprint.} 
%
%
 \vspace*{-.1in}
\item Bloch, A., N. Leonard and J. Marsden (1998).
Matching and stabilization by the method of controlled 
Lagrangians. {\it Proc. IEEE Conf. on Decision and Control},
Tampa, FL, pp. 1446-1451.
\vspace*{-.1in}
\item Bloch, A., N. Leonard and J. Marsden (1999).
 Stabilization of the pendulum on a rotor arm 
by the method of controlled Lagrangians. 
{\it Proc. IEEE Int. Conf. on Robotics and Automation}, 
 Detroit, MI, pp. 500-505.
%
\vspace*{-.1in}
\item  Hamberg, J.(1999). General matching conditions in the theory of 
controlled Lagrangians. 
{\it Proceedings of 
the 38th Conference on Decision and Control},
 Phoenix, 
 AZ.
 \vspace*{-.1in}
\item van der Schaft, A. J. (1986). 
 Stabilization of Hamiltonian systems.
 {\it Nonlinear Analysis, Theory, Methods \& Applications}, 
{\bf 10},   1021-1035.
 
\end{description}

\end{document}